\documentclass[12pt]{amsart}
\theoremstyle{plain}
\newtheorem{Thm}{Theorem}
\newtheorem{Lem}[Thm]{Lemma}
\newtheorem{ass}[Thm]{Assertion}

\errorcontextlines=0 \numberwithin{equation}{section}

\begin{document}
\large
\title[$L^2$ Forms and Ricci flow ]
{$L^2$ Forms and Ricci flow with bounded curvature on Complete
Non-compact manifolds }
\author{Li MA and Yang Yang}

\address{Department of mathematical sciences \\
Tsinghua university \\
Beijing 100084 \\
China}

\email{lma@math.tsinghua.edu.cn} \dedicatory{}
\date{Oct. 15th, 2002}

\keywords{Ricci flow, forms, maximum principle.} \subjclass{53C.}

\begin{abstract}
In this paper, we study the evolution of $L^2$ one forms under Ricci
flow  with bounded curvature on a non-compact Rimennian manifold. We
show on such a manifold that the $L^2$ norm of a smooth one form
with compact support is non-increasing along the Ricci flow with
bounded curvature. The $L^{\infty}$ norm is showed to have
monotonicity property too. Then we use $L^{\infty}$ cohomology of
one forms with compact support to study the singularity model for
the Ricci flow on $S^1\times \mathbb{R}^{n-1}$.
\end{abstract}
\maketitle

\section{ Introduction}
In \cite{P}, G.Perelman found remarkable functionals, which are
monotone along Ricci flow, to study the singularity model of Ricci
flow on compact Riemannian manifolds. See also \cite{KL} and
\cite{CM05} for beautiful illustration and and applications. The
story will be different for complete non-compact Riemannian
manifolds. In this paper, we try to study the singularity model for
the Ricci flow on a non-compact Riemannian manifold $(M,g_0)$ with
bounded curvature. So we consider monotonicity properties related to
Ricci flow and other heat equations. Interestingly, L.Ni \cite{N04}
find a very nice monotonicity formula for heat flow based on the
fundamental works of G.Perelman \cite{P} and Li-Yau \cite{LY}. The
basic tool in our mind is to use invariants of Riemannian manifolds,
such as that related to harmonic forms, so we shall consider the
behavior of an $L^2$ co-homology class of one forms, saying $\Phi$,
under the Ricci flow. For background of $L^2$ harmonic forms, one
may refer to \cite{Do}, \cite{An}, \cite{CG}, and \cite{G91}. Here,
we just say that the $L^2$ space of one forms is the completion of
the space of smooth one forms with compact support with respect to
the $L^2$ inner product on $(M,g)$. M.Gaffney \cite{Ga54} proved
that the Hodge decomposition theorem is still true on complete
Riemannian manifolds. We mention another interesting result that
$L^2$ harmonic functions are constants on complete Riemannian
manifolds \cite{Y76}.

By definition, the Ricci-Hamilton flow is
  the evolution equation for Riemannian metrics:
  $$
\partial_tg_{ij}=-2R_{ij},\;\;\;\mbox{on $M_T:=M\times [0,T)$}
  $$
where $R_{ij}$ is the Ricci tensor of the metric $g:=g(t)=(g_{ij})$
in local coordinates $(x^i)$ and $T$ is the maximal existing time
for the flow. The existence of Ricci flow with bounded sectional
curvature on a complete non-compact Riemannian manifold had been
established by Shi \cite{Shi89} in 1989, and this is a very useful
result in
 Riemannian Geometry. Interestingly, the maximum principle of heat
 equation is true on such a flow. Then we can easily show that
 the Ricci flow
preserves the property of nonnegative scalar curvature (see also
\cite{H95}). Given a smooth $L^2$ one form $\phi$ with compact
support on a Riemannian manifold $(M,g)$. Recall that its norm is
defined by
$$
||\phi||_{L^2g}=(\int_M|\phi|^2_g(x)dv_g)^{1/2}.
$$
Assume that $d\phi=0$. Let $\Phi=[\phi]$ be the $L^2$ cohomology
class of the form $\phi$ in $(M,g)$. Define
$$
||\Phi||_{L^2g}=inf_{\varphi\in \Phi}||\varphi||_{L^2g}
$$
and
$$||\Phi||_2(t)=||\Phi||_{L^2g(t)}$$
for the flow $\{g(t)\}$.  It is well-known that $||\Phi||_{L^2g}$ is
a norm on $H^1_{dR}(M,\mathbf{R})$. We denote by $d_g(x,y)$ the
distance of two points $x$ and $y$ in $(M,g)$.

Our first result is

\begin{Thm}\label{yy} Let $(M,g_0)$ be a complete
non-compact Riemannian manifold
with non-negative scalar curvature. Assume that $g(t)$ is a Ricci
flow with bounded curvature on $[0,T)$ with initial metric
$g(0)=g_0$ on $M$. Then for a $L^2$ one form $\phi$, we have
$$
||\phi||_{L^2g(t)}\leq ||\phi||_{L^2g(s)},\;\;{for}\;\;t>s,
$$
along the Ricci flow $g(t)$.
\end{Thm}

 Our second result is the following
\begin{Thm}\label{ml}
Assume $\alpha\in H_1(M,\mathbf{Z})$ is an element of infinite
order. Assume that $T>0$ and  $g(t)$ is a Ricci flow with bounded
curvature on $[0,T)$ with initial metric $g(0)=g_0$ on a non-compact
complete Riemannian manifold $(M,g_0)$. Then there is a uniform
constant $c=(\alpha,g(0),T)$ such that
$$ L_{\alpha}(g(t))\geq c
$$
for all $t\in [0,T)$. Here $L_{\alpha}(g(t))$ is the infimum of the
length of the curves $a\in\alpha$ in the metric $g(t)$.
\end{Thm}

One example with such an $\alpha$ is $M=S^1\times
 \mathbf{R}^{n-1}$. For the compact case, this is an interesting
 result of Ilmanen and Knopf\cite{IK}.

\textbf{{Definition.}} \emph{We say that the Riemannian manifold
$(M_{\infty},h)$ is the infinite time  blow up limit of Ricci flow
$(M,g(t))$ if there exist a sequence $t_k\to \infty$, a sequence
$\lambda_k\to\infty$, and a point sequence $(x_k)\subset M$ with
limit $x\in M_{\infty}$ with the blow up metric sequence $g_k(t)$,
which is defined by
$$
g_k(t)=\lambda_kg(t_k+\frac{t}{\lambda_k}),
$$
such that $(M,g_k(t),x_k)$ converges locally smoothly to the limit
$(M_{\infty},h,x)$. }

Note that for a smooth closed curve $\Gamma$, the lengths of
$\Gamma$ in the metric $g(t_k)$ and $g_k(0)$ related by
$$L(\Gamma,g_k(0))=\lambda_kL(\Gamma,g(t_k)).$$
In particular, if $L(\Gamma,g(t))\geq c>0$, one has
$$
L(\Gamma,g_k(\cdot))\geq \lambda c\to+\infty.
$$

 Using above and Theorem \ref{ml}, we know (as in \cite{IK}) that
\begin{ass}
 $M=S^1\times \mathbf{R}^{n-1}$ with metric $\bar{g}(t)$
 can not be the infinite time blow up
 limit of Ricci flow with \emph{curvature decay at infinity}. Here
 $$\bar{g}(t)=ds^2+g_{can}$$
 with $g_{can}$ being the standard metric on $\mathbf{R}^{n-1}$.
\end{ass}
The \emph{curvature decay at infinity} for Ricci flow means that
each Riemannian metric $g(t)$ has curvature decay of positive order
at infinity. Recall that, for bounded smooth tensor $\phi$ on a
complete $(M,g)$, the meaning of \emph{the decay (of order
$\sigma>0$) at infinity} is that
 for
$d_g(x,o)\to+\infty$, we have
$$d_g(x,0)^{\sigma}|\phi|_{g(x)}\to 0.$$

The Ricci flow preserves the curvature decay at infinity, which was
proved by Dai and Ma in \cite{DM1}. Since along the Ricci flow, we
have that the decay at infinity is preserved, and the proof of above
result is basically reduced to compact domains. So the proof for
compact case in \cite{IK} can be carried to our case, so the detail
shall be omitted.

We introduce some notations. Let $\phi=\phi_jdx^j$ in local
coordinates $(x^i)$. Given Riemannian metric $g$, we denote by
$d_g(x,o)$ be the distance between two points $x$ and $o$. Set
$(g^{ij})=(g_{ij})^{-1}$ and $R^j_i=g^{jk}g_{ki}$.

\section{Proof of Theorem \ref{yy}}
In this section, we first prove a useful lemma.
\begin{Lem} \label{zhao} The $L^1$ property of
 non-negative bounded solutions to the scalar differential inequality
$$
u_t\leq\Delta u
$$
is preserved along the Ricci flow with bounded curvature and
non-negative scalar curvature. Furthermore,
$$
u(t)\leq u(s),\;\;{for}\;\;t>s.
$$
\end{Lem}
\begin{proof}
Recall that on $M$,
$$
u_t\leq\Delta u.
$$
Let $p=1+\epsilon$ with small $\epsilon>0$. Take $r>0$ large. Let
$\eta$ be a non-negative cut-off function such that $0\leq \eta\leq
1$ on $M$, $\eta=1$ on $B_r(o)$, $\eta=0$ outside $B_{2r}(o)$, and
$$
|\nabla \eta|^2\leq 4\eta/r^2.
$$
Then
\begin{align*}
&\int_0^Tdt\int_M\eta^2u^{p-1}(u_t)\leq
\int_0^Tdt\int_M\eta^2u^{p-1}\Delta u\\
&=-2\int_0^Tdt\int_M\eta
u^{p-1}<\nabla\eta,\nabla u>\\
&-(p-1)\int_0^Tdt\int_M\eta^2u^{p-2}|\nabla u|^2\\
&\leq\frac{2}{p-1}\int_0^Tdt\int_M|\nabla
\eta|^2u^p-\frac{2(p-1)}{p^2}\int_0^Tdt\int_M\eta^2|\nabla(u^{p/2})|^2\\
&\leq \frac{2}{(p-1)r^2}\int_0^Tdt\int_M\eta u^p
-\frac{2(p-1)}{p^2}\int_0^Tdt\int_M\eta^2|\nabla(u^{p/2})|^2\\
&\leq \frac{2}{(p-1)r^2}\int_0^Tdt\int_M\eta u^p\to 0
\end{align*}
as $r\to \infty$. Here we have used the decay condition that $u\in
L^p$ for every $p>1$ small.

 By a direct computation, we have
\begin{align*}
&\int_0^Tdt\int_M\eta^2u^{p-1}u_t\\
&=\frac{1}{p}\int_0^Tdt\frac{d}{dt}\int_M
\eta^2u^p\\
&+\frac{1}{p}\int_0^Tdt\int_M \eta^2u^{p}R\\
&=\frac{1}{p}\int_M \eta^2u^p(t)
-\frac{1}{p}\int_M \eta^2u^p(0)\\
&+\frac{1}{p}\int_0^Tdt\int_M \eta^2u^{p}R.
\end{align*}
Hence, we have
\begin{align*}
&\frac{1}{p}\int_M\eta^2u^p(t)\\
&-\frac{1}{p}\int_M \eta^2u^p(0) +\frac{1}{p}\int_0^Tdt\int_M
\eta^2u^{p}R\leq o(1)\to 0.
\end{align*}
 Sending $r\to +\infty$, we get that
\begin{align*}
\frac{1}{p}\int_M u^p(t) -\frac{1}{p}\int_M
u^p(0)+\frac{1}{p}\int_0^Tdt\int_M u^{p}R\leq 0.
\end{align*}
Sending $p\to 1$, we have that
\begin{align*}
\int_M u(t) -\int_M u(0)+\int_0^Tdt\int_M u^{1}R\leq 0.
\end{align*}
That is,
$$
\int_M u(t) -\int_M u(0)+\int_0^Tdt\int_M u^{1}R\leq 0,
$$
which implies that $u\in L^1$ for each $t>0$. The monotonicity of
$u$ comes clearly from the argument above.
\end{proof}

We remark that for a compact Riemannian manifold, the corresponding
result of above lemma is also true.

Let $\phi_0\in \Phi$ be a closed one-form with compact support. We
try to evolve $\phi_0$ by the family of 1-forms $\phi(t)$ satisfying
the heat equation
$$
\frac{\partial}{\partial t}\phi=\Delta_d\phi
$$
with $\phi(t)=\phi_0$ at $t=0$. Here $-\Delta_d=\delta_c
d+d_c\delta$ is the Hodge-DeRham Laplacian of $g(t)$ in the sense of
M.P.Gaffney \cite{Ga54} (see also \cite{SY94}). We now recall the
famous Bochner-formula for 1-form.  Then we have
$$
\Delta_d\phi_k=\Delta \phi_k-R^j_k\phi_j,
$$
where $\Delta= \nabla^*\nabla$ is the rough Laplacian of $g(t)$.
Clearly, we have
\begin{align}
\frac{1}{2}\frac{\partial}{\partial
t}|\phi|^2&=\frac{1}{2}\frac{\partial}{\partial
t}g^{ij}\phi_i\phi_j\nonumber\\
&=R^{ij}\phi_{i}\phi_{j}\phi^i+\phi^j\frac{\partial}{\partial
t}\phi_j\nonumber\\
&=\phi^j\Delta \phi_j \nonumber\\
&=\frac{1}{2}\Delta |\phi|^2-|\nabla\phi|^2.\label{max}
\end{align}

It is clear that Lemma \ref{zhao} implies Theorem \ref{yy}.

Then we we have the following consequence of Theorem \ref{yy}
\begin{ass}
 Let $(M,g_0)$ be a complete non-compact Riemannian manifold
with non-negative scalar curvature. Assume that $g(t)$ is a Ricci
flow with bounded curvature on $[0,T)$ with initial metric
$g(0)=g_0$ on $M$. Then for an $L^2$ harmonic one form $\phi$, we
have
$$
||[\phi]||_{L^2g(t)}\leq ||[\phi]||_{L^2g(s)},
\;\;{for}\;\;t>s\geq0,
$$
along the Ricci flow $g(t)$.
\end{ass}

The proof is just the combination of definition of norms of $[\phi]$
and Theorem \ref{yy}.

\section{proof of Theorem \ref{ml}}
Let $\Phi=[\phi]$ on $(M,g)$ as in last section.
 Define, for $\varphi\in \Phi$,
$$
||\varphi||_{g}=\sup_{x\in M}|\varphi(x)|_g
$$
 and
$$
|\Phi|_{g}=\inf_{\varphi\in\Phi}||\varphi||_{g}.
$$
We write by
$$
|\Phi|(t)=|\Phi|_{g(t)},
$$
along the flow $\{g(t)\}$.

We now consider the equation \ref{max} again. As in \cite{DM1} (see
also Theorem 18.2 in Hamilton \cite{H95}), using the maximum
principle to \ref{max} we know that $|\phi(t)|_{g(t)}$ has a decay
of any order at infinity and
$$
||\phi(t)||_{g(t)}\leq ||\phi(0)||_{g(0)},\label{mali}
$$
for $t\in [0,T)$.

Then we claim that

\begin{ass}
 $|\Phi|_g$  is a
norm on $H^1_{dR}(M,\mathbf{R})$.
\end{ass}
\begin{proof}
The triangle inequalities and homogeneity are straight forward. We
only need to show the positivity. Assume that  $|\Phi|_g=0$. Then we
can choose a sequence of smooth forms $\{\phi_k\}$ with compact
support in the class $\Phi$ such that $||\phi_k||_{g}\to 0$. Take
any smooth $\phi\in \Phi$ with compact support. So we can assume
that $$ ||\phi_k||_{g}\leq C
$$
for some uniform constant. Then for each $k$, we can write by
$$ \phi-\phi_k=dF_k
$$
for some smooth function $F_k$ with compact support (see
\cite{Y76}). Hence
$$
|dF_k|\leq C+||\phi ||_{g}.
$$
By adding a constant to each $F_k$, we can assume that $F_k(o)=0$.
Then we can have an uniformly locally convergent subsequence, still
denoted by $\{F_k\}$, with limit $F$ a Lipschitz function with
$F(o)=0$. This implies that
\begin{align*}
ess \sup|\phi-dF|_g&\leq \sup\lim_{k\to +\infty}||\phi-dF_k||_{g}\\
&=\sup\lim_{k\to +\infty}||\phi_k||_{g}=0.
\end{align*}
Then we conclude that $\phi=dF$ almost everywhere. Since $\phi$ is
smooth, we know that $F$ is also smooth. Hence $\Phi=[\phi]=0$.
\end{proof}

Along the Ricci flow, we also have the following monotonicity
formulae.

\begin{ass}
Under the Ricci flow with bounded curvature, $|\Phi|(t)$ is
non-increasing.
\end{ass}
\begin{proof} We only prove the monotonicity for $|\Phi|(t)$.
the other is the same. For every $\epsilon>0$, we choose a smooth
one form $\phi_0\in \Phi$ with compact support such that
$$
||\phi_0||_{g_0}\leq |\Phi|_{g_0}+\epsilon.
$$
We now study the equation
\begin{equation}
\frac{\partial}{\partial t}\phi=\Delta_d\phi \label{yanyang}
\end{equation}
with initial data $\phi(0)=\phi_0$. As we show in \cite{DM1}, we
know that the solution $\phi(t)$ exists as long as $g(t)$ exists,
and the solution is uniformly bounded. It is easy to verify that
$$
\phi(t)\in \Phi.
$$
In fact, we can solve the following diffusion equation
$$
\partial_t F=\Delta F-\delta\phi_0,\;\;F(0)=0
$$
and get the unique smooth solution $F(t)$ with decay (in $L^2$).
Note that $\phi_0+dF$ solves (\ref{yanyang}) with the same initial
data $\phi_0$. By the maximum principle we conclude that
$\phi=\phi_0+dF$ on $M$. Hence, by (\ref{mali}) we have
$$
|\Phi|(t)\leq ||\phi||_{g(t)}\leq ||\phi_0||_{g(0)}\leq
|\Phi|_{g_0}+\epsilon.
$$
\end{proof}

We now give the proof of Theorem \ref{ml}.

\begin{proof}
As in \cite{IK}, we can choose a closed curve $a$ representing
$\alpha$ and a smooth one form $\phi\in \Phi$ with compact support.
Then
$$
0<<\Phi,\alpha>=\int_{a}\phi\leq ||\phi||_{g(t)}\cdot L(a,g(t)).
$$
Here $L(a,g(t))$ is the length of the curve $a$ in the metric
$g(t)$. Since $\phi$ is arbitrary, we have
$$
0<<\Phi,\alpha>\leq |\Phi|(t)L(a,g(t)).
$$

Using the monotonicity of $|\Phi|(t)\leq |\Phi|(0)$, we get that
$$
L(a,g(t))\geq \frac{<\Phi,\alpha>}{|\Phi|(0)}>0.
$$
\end{proof}

\end{document}